\documentclass[10pt]{amsart}

\usepackage{amsmath,amsthm,amsfonts,amscd,amssymb}
\usepackage{hyperref,wasysym}

\newtheorem{thm}{Theorem}[section]
\newtheorem{cor}[thm]{Corollary}
\newtheorem{lem}[thm]{Lemma}
\newtheorem{prop}[thm]{Proposition}

\newtheorem{ex}{Example}

\theoremstyle{definition}

\theoremstyle{remark}
\newtheorem{rem}{Remark}[section]

\begin{document}

\title{On the geometry of nearly orthogonal lattices}
\author{Lenny Fukshansky}
\author{David Kogan}\thanks{Fukshansky was partially supported by the Simons Foundation grant \#519058}

\address{Department of Mathematics, 850 Columbia Avenue, Claremont McKenna College, Claremont, CA 91711}
\email{lenny@cmc.edu}
\address{Institute of Mathematical Sciences, Claremont Graduate University, Claremont, CA 91711}
\email{david.kogan@cgu.edu}


\begin{abstract}
Nearly orthogonal lattices were formally defined in~\cite{baraniuk}, where their applications to image compression were also discussed. The idea of ``near orthogonality" in $2$-dimensions goes back to the work of Gauss. In this paper, we focus on well-rounded nearly orthogonal lattices in~$\mathbb R^n$ and investigate their geometric and optimization properties. Specifically, we prove that the sphere packing density function on the space of well-rounded lattices in dimension $n\geq 3$ does not have any local maxima on the nearly orthogonal set and has only one local minimum there: at the integer lattice~$\mathbb Z^n$. Further, we show that the nearly orthogonal set cannot contain any perfect lattices for~$n \geq 3$, although it contains multiple eutactic (and even strongly eutactic) lattices in every dimension. This implies that eutactic lattices, while always critical points of the packing density function, are not necessarily local maxima or minima even among the well-rounded lattices. We also prove that a (weakly) nearly orthogonal lattice in~$\mathbb R^n$ contains no more than~$4n-2$ minimal vectors (with any smaller even number possible) and establish some bounds on coherence of these lattices. 
\end{abstract}

\maketitle

\maketitle

\def\A{{\mathcal A}}
\def\B{{\mathcal B}}
\def\C{{\mathcal C}}
\def\D{{\mathcal D}}
\def\F{{\mathcal F}}
\def\x{{\mathcal H}}
\def\I{{\mathcal I}}
\def\J{{\mathcal J}}
\def\K{{\mathcal K}}
\def\L{{\mathcal L}}
\def\M{{\mathcal M}}
\def\O{{\mathcal O}}
\def\R{{\mathcal R}}
\def\s{{\mathcal S}}
\def\V{{\mathcal V}}
\def\W{{\mathcal W}}
\def\X{{\mathcal X}}
\def\Y{{\mathcal Y}}
\def\H{{\mathcal H}}
\def\OO{{\mathcal O}}
\def\BB{{\mathbb B}}
\def\cee{{\mathbb C}}
\def\pee{{\mathbb P}}
\def\que{{\mathbb Q}}
\def\real{{\mathbb R}}
\def\zed{{\mathbb Z}}
\def\hyp{{\mathbb H}}
\def\aa{{\mathfrak a}}
\def\qbar{{\overline{\mathbb Q}}}
\def\eps{{\varepsilon}}
\def\ahat{{\hat \alpha}}
\def\bhat{{\hat \beta}}
\def\gt{{\tilde \gamma}}
\def\h{{\tfrac12}}
\def\be{{\boldsymbol e}}
\def\bei{{\boldsymbol e_i}}
\def\bff{{\boldsymbol f}}
\def\ba{{\boldsymbol a}}
\def\bb{{\boldsymbol b}}
\def\bc{{\boldsymbol c}}
\def\bm{{\boldsymbol m}}
\def\bk{{\boldsymbol k}}
\def\bi{{\boldsymbol i}}
\def\bl{{\boldsymbol l}}
\def\bq{{\boldsymbol q}}
\def\bu{{\boldsymbol u}}
\def\bt{{\boldsymbol t}}
\def\bs{{\boldsymbol s}}
\def\bv{{\boldsymbol v}}
\def\bw{{\boldsymbol w}}
\def\bx{{\boldsymbol x}}
\def\bX{{\boldsymbol X}}
\def\bz{{\boldsymbol z}}
\def\bwy{{\boldsymbol y}}
\def\bY{{\boldsymbol Y}}
\def\bL{{\boldsymbol L}}
\def\baa{{\boldsymbol\alpha}}
\def\bbb{{\boldsymbol\beta}}
\def\bet{{\boldsymbol\eta}}
\def\bxi{{\boldsymbol\xi}}
\def\bo{{\boldsymbol 0}}
\def\bol{{\boldkey 1}_L}
\def\ep{\varepsilon}
\def\p{\boldsymbol\varphi}
\def\q{\boldsymbol\psi}
\def\rank{\operatorname{rank}}
\def\aut{\operatorname{Aut}}
\def\lcm{\operatorname{lcm}}
\def\sgn{\operatorname{sgn}}
\def\spn{\operatorname{span}}
\def\md{\operatorname{mod}}
\def\Norm{\operatorname{Norm}}
\def\dim{\operatorname{dim}}
\def\det{\operatorname{det}}
\def\Vol{\operatorname{Vol}}
\def\rk{\operatorname{rk}}
\def\Gal{\operatorname{Gal}}
\def\WR{\operatorname{WR}}
\def\WO{\operatorname{WO}}
\def\GL{\operatorname{GL}}

\section{Introduction}
\label{intro}

Let $L$ be a lattice of rank $n \geq 2$ in $\real^n$, and define its {\it minimal norm} to be
$$|L| = \min \left\{ \|\bx\| : \bx \in L \setminus \{\bo\} \right\},$$
where $\|\ \|$ is the Euclidean norm on $\real^n$. Then its set of {\it minimal vectors} is
$$S(L) = \left\{ \bx \in L : \|\bx\| = |L| \right\}$$
We define {\it coherence} of the lattice $L$ to be
$$C(L) = \max \left\{ |\cos \aa(\bx,\bwy)| : \bx, \bwy \in S(L), \bx \neq \pm \bwy \right\},$$
where we write $\aa(\bx,\bwy)$ for the angle between the two vectors~$\bx$ and~$\bwy$. Notice that
$$0 \leq C(L) \leq 1/2,$$
since the angle between any two minimal vectors $\bx \neq \pm \bwy$ of a lattice is in the interval $[\pi/3,2\pi/3]$ (see, for instance, Lemma~3.1 of~\cite{hex}). The coherence of lattices was recently introduced in~\cite{ejc} by analogy with coherence of frames, an important notion in signal processing (see, for instance, \cite{bajwa}, \cite{mixon}, and the references therein). Generally, one can think of coherence of a set of vectors as a measure of how far are they from being orthogonal. Various techniques in error-correcting codes and signal recovery employ overdetermined equal-norm spanning sets in Euclidean vector spaces that are as close to orthogonal as possible. Since lattices are frequently used in signal processing and digital communications, we want to better understand the coherence properties of such sets coming from the minimal vectors of lattices. 

A classical optimization problem studied on lattices is the sphere packing problem. Let us write $\det(L)$ for the determinant of the lattice $L$: it is computed as the absolute value of the determinant of a basis matrix of $L$ and is equal to the volume of the corresponding fundamental parallelotope spanned by the vectors of this basis. This definition does not depend on the choice of the basis, since change of basis corresponds to multiplication by a matrix from $\GL_n(\zed)$. There is a sphere packing associated with every lattice $L$: it consists of spheres of radius $|L|/2$ centered at the points of $L$, and its {\it density} is the proportion of the space it occupies, which can be computed as
$$\delta(L) := \frac{\omega_n |L|^n}{2^n \det L},$$
where $\omega_n$ is the volume of a unit ball in~$\real^n$ (see \textsection{1.8} of \cite{martinet}). The space of lattices in $\real^n$ can be identified with $\GL_n(\real) / \GL_n(\zed)$, i.e., nonsingular real $n \times n$ matrices modulo right multiplication by nonsingular integer $n \times n$ matrices, and $\delta$ is a continuous function on this space. Lattices that are local maxima of $\delta$ are called {\it extreme}. In~\cite{ejc} some heuristics were presented, speculating that there may be an inverse correlation between the coherence and packing density on lattices. Our goal here is to investigate this correlation on the important class of nearly orthogonal lattices.

Two lattices $L_1,L_2$ in $\real^n$ are called {\it similar} (written $L_1 \sim L_2$) if there exists a positive constant $\alpha$ and an $n \times n$ real orthogonal matrix $U$ such that $L_2 = \alpha U L_1$. This is an equivalence relation on the space of lattices in $\real^n$ with equivalence classes referred to as {\it similarity classes} of lattices in~$\real^n$. The space of similarity classes of lattices in $\real^n$ can be identified with $\left( \real_{+} \times \O_n(\real) \right) \backslash \GL_n(\real) / \GL_n(\zed)$, i.e., lattices modulo left multiplication by positive constants and orthogonal matrices. Metric topology on this space is induced by the usual Euclidean norm on the space of $n \times n$ real matrices viewed as vectors in~$\real^{n^2}$: it is the quotient topology on the space of lattices $\GL_n(\real) / \GL_n(\zed)$ modulo the left-multiplication action of $\real_{+} \times \O_n(\real)$, under which the similar lattices are identified, and hence correspond to the same point (see Chapter~1 of \cite{martinet} for a more detailed description of topology on the space of lattices and the space of similarity classes of lattices, especially \textsection{1.1} and \textsection{1.7}). It is easy to notice that if $L_1 \sim L_2$, then $\delta(L_1) = \delta(L_2)$ and $C(L_1) = C(L_2)$. In particular, extreme lattices can only be similar to extreme lattices, and $\delta$ is a continuous function on the space of similarity classes of lattices in $\real^n$.

A lattice $L$ is called {\it well-rounded} (abbreviated WR) if 
$$\spn_{\real} S(L) = \spn_{\real} L.$$
Because WR lattices can only be similar to WR lattices, we write $\WR_n$ for the space of similarity classes of WR lattices. It is a well-known fact that extreme lattices must be WR, which is why the study of the packing density function on the space of lattices is usually restricted to WR lattices. Detailed further information on lattices, their geometric properties, and packing density can be found in the classical texts~\cite{conway:sloane} and~\cite{martinet}.

We focus on the important class of nearly orthogonal lattices as defined in~\cite{baraniuk}. These lattices appear to be useful in image processing, signal recovery, and related areas (see, for instance~\cite{baraniuk} and~\cite{chen}). Let $B = \{ \bb_1,\dots,\bb_n \}$ be an ordered basis for a lattice $L$ in~$\real^n$, and define a sequence of angles $\theta_1,\dots,\theta_{n-1}$ as follows: each $\theta_i$ is the angle between $\bb_{i+1}$ and the subspace $\spn_{\real} \{ \bb_1,\dots,\bb_i \}$. It is then clear that each $\theta_i \in (0,\pi/2]$. For a given value $\theta \in [0,\pi/2]$, we will say that $B$ is a {\it weakly $\theta$-orthogonal} basis if $\theta_i \geq \theta$ for each $1 \leq i \leq n-1$. A basis $B$ is called {\it $\theta$-orthogonal} if every ordering of it is weakly $\theta$-orthogonal. Notice that if some lattice $L$ has a (weakly) $\theta$-orthogonal basis, then so does every lattice in its similarity class: we will call such lattices (weakly) $\theta$-orthogonal. Let us then write $\WO_n(\theta)$ for the space of similarity classes of all weakly $\theta$-orthogonal lattices, and $\WO^*_n(\theta)$ for the space of similarity classes of all $\theta$-orthogonal lattices in~$\real^n$. We will also write simply $\WO_n$ (respectively, $\WO^*_n$) for $\WO_n(\pi/3)$ (respectively, $\WO^*_n(\pi/3)$), which will be especially important to us; we will call lattices in~$\WO_n$ (respectively, $\WO^*_n$) {\it weakly nearly orthogonal} (respectively, {\it nearly orthogonal}) and refer to their corresponding (weakly) $\pi/3$-orthogonal bases as (weakly) nearly orthogonal bases. It is shown in~\cite{baraniuk} that if $L$ has a weakly nearly orthogonal basis $B$, then $B$ contains a minimal vector of $L$. As discussed in~\cite{baraniuk}, not every lattice has a weakly nearly orthogonal basis. The authors of~\cite{baraniuk} established certain specific properties of such bases, when they exist (generalizing a classical 2-dimensional result of Gauss) with a view towards applications in image compression. It is our goal to move further and investigate the geometric and optimization properties of WR lattices possessing such bases. Let us define 
$$\W_n(\theta) = \WR_n \cap \WO_n(\theta) \text{ and } \W^*_n(\theta) = \WR_n \cap \WO^*_n(\theta),$$
i.e., the set of similarity classes of WR lattices in~$\real^n$ that have a (weakly) $\theta$-orthogonal basis; we will write simply $\W_n$ for $\W_n(\pi/3)$ (respectively, $\W^*_n$ for $\W^*_n(\pi/3)$). 

Each similarity class can be represented by a lattice with minimal norm~$1$: from now on, we will always use such representatives.  When we write $L \in \W_n(\theta)$, we will mean the similarity class of $L$ where $|L|=1$. For a given basis $B$ of a lattice $L$, we define
\begin{equation}
\label{mu_nu}
\begin{aligned}
& \mu(B) := \min \{ | \cos \aa(\bb_i, \bb_j) | : 1 \leq i \neq j \leq n \},\\
& \nu(B) := \max \{ | \cos \aa(\bb_i, \bb_j) | : 1 \leq i \neq j \leq n \}).
\end{aligned}
\end{equation}
In this paper, we investigate geometric properties and packing density of WR nearly orthogonal lattices. To do so, we examine the minimal vectors of these lattices and the angles between them. Here is our first result in this direction.

\begin{thm} \label{dense_coh} Let $L \in \W^*_n$ with a nearly orthogonal basis $B = \{ \bb_1,\dots,\bb_n \}$. Then $B \subseteq S(L)$, and so $\mu(B) \leq \nu(B) \leq C(L)$. Let $\eps > 0$ and let $\BB_{\eps}(L)$ be a ball of radius $\eps$ centered at $L$ in the space of similarity classes of lattices in~$\real^n$. The following statements are true:

\begin{enumerate}

\item If $\mu(B) < 1/2$, then there exists $L' \in \W^*_n \cap \BB_{\eps}(L)$ with nearly orthogonal basis $B'$ such that $\mu(B') > \mu(B)$ and
$$\delta(L') = \frac{\sqrt{1 - \mu(B)^2}}{\sqrt{1 - \mu(B')^2}}\ \delta(L) > \delta(L).$$

\item If $\nu(B) > 0$, then there exists $L'' \in \W^*_n \cap \BB_{\eps}(L)$ with nearly orthogonal basis $B''$ such that $\nu(B'') < \nu(B)$ and
$$\delta(L'') = \frac{\sqrt{1 - \nu(B)^2}}{\sqrt{1 - \nu(B'')^2}}\ \delta(L) < \delta(L).$$

\item If $n \geq 3$, then $\mu(B) < 1/2$ for every nearly orthogonal basis $B$ of any lattice $L \in \W^*_n$. If $n=2$, $\mu(B) = 1/2$ if and only if $L$ is the hexagonal lattice, in which case this is true for any nearly orthogonal basis $B$ of $L$.

\item $C(L) = 0$ if and only if $L = \zed^n$.

\end{enumerate}
\end{thm}

\noindent
Loosely speaking, Theorem~\ref{dense_coh} asserts that a lattice $L \in \W^*_n$ can be locally modified to increase or decrease the packing density. We prove Theorem~\ref{dense_coh} in Section~\ref{dense}. In fact, as becomes evident from our proof, such local modification simultaneously increases or decreases both, packing density and coherence of the lattice. Here is an immediate consequence of this theorem.

\begin{cor} \label{extreme} If $n \geq 3$, then $\W^*_n$ does not contain any extreme lattices. If $n = 2$, the hexagonal lattice is the unique extreme lattice, which is contained in~$\W^*_2$. On the other hand, $\W^*_n$ for every $n \geq 2$ contains a unique minimum of the packing density function on the set of WR lattices, the integer lattice~$\zed^n$.
\end{cor}

The fact that $\W_n(\theta)$ with $\theta > \pi/3$ cannot contain extreme lattices already follows from results of~\cite{baraniuk} in a different manner. A lattice $L$ is called {\it weakly eutactic} if there exist real numbers $c_1,\dots,c_n$, called {\it eutaxy coefficients}, such that
\begin{equation}
\label{eutaxy}
\|\bv\|^2 = \sum_{\bx \in S(L)} c_i  (\bv,\bx_i)^2
\end{equation}
for all $\bv \in \real^n$, where $(\ ,\ )$ stands for the usual scalar product of vectors. If eutaxy coefficients are positive, $L$ is called {\it eutactic}, and if $c_1 = \dots = c_n > 0$, the lattice~$L$ is called {\it strongly eutactic}; for instance, the integer lattice~$\zed^n$ is strongly eutactic. Further, $L$ is {\it perfect} if the set $\{ \bx \bx^{\top} : \bx \in S(L) \}$ spans the space of $n \times n$ real symmetric matrices; here and throughout vectors are always written as columns. Both, eutactic and perfect lattices are necessarily WR, and eutaxy and perfection properties are preserved on similarity classes, same as well-roundedness. A famous theorem of Voronoi (1908, \cite{voronoi}) asserts that $L$ is extreme if and only if $L$ is eutactic and perfect. Notice that in order for $L$ to be perfect $S(L)$ needs to contain at least $\frac{n(n+1)}{2}$ pairs of $\pm$ minimal vectors, the dimension of the space of $n \times n$ real symmetric matrices. On the other hand, if $L \in \W_n(\theta)$ with $\theta > \pi/3$ and $B$ is its weakly $\theta$-orthogonal basis, then Corollary~1 of~\cite{baraniuk} states that $S(L) = \pm B$. Hence $L$ cannot be perfect, since $n < \frac{n(n+1)}{2}$ for all $n \geq 2$. This, however, does not imply our result for~$\W^*_n$. Indeed, while we do prove that $B \subseteq S(L)$ for any nearly orthogonal basis $B$ of any lattice $L \in \W^*_n$, it is possible to construct lattices in~$\W^*_n$ with larger sets of minimal vectors.

\begin{thm} \label{SL_more} Let $n \geq 2$. For each $0 \leq m \leq [n/2]$ there exists a strongly eutactic lattice $L_{n,m} \in \W^*_n$ with
$$|S(L_{n,m})| = 2(n+m).$$
In particular, if $m=[n/2]$, then
$$|S(L_{n,m})| = \left\{ \begin{array}{ll}
3n & \mbox{if $n$ is even} \\
3n-1 & \mbox{if $n$ is odd}.
\end{array}
\right.$$
Furthermore, there exist $L \in \W_n$ such that $|S(L)|$ is any even number between $3n$ and $4n-2$, inclusive. On the other hand, $|S(L)| \leq 4n-2$ for every~$L \in \W_n$, and if $|S(L)| > 3n$ then $L$ cannot be in~$\W^*_n$.
\end{thm}

\noindent
We prove Theorem~\ref{SL_more} in Section~\ref{SL}, in particular constructing explicit families of lattices. Notice that the set $\W^*_n$ contains strongly eutactic lattices. Further, in Example~\ref{eutactic_3d} we exhibit a~$3$-dimensional irreducible non-perfect eutactic (but not strongly eutactic) lattice, which is in~$\W_3$, but not in~$\W^*_3$; recall that a nonzero lattice is {\it irreducible} if it cannot be represented as an orthogonal direct sum of two proper sublattices. On the other hand, since perfect lattices must have at least~$n(n+1)$ minimal vectors, Theorem~\ref{SL_more} implies an immediate corollary.

\begin{cor} \label{no_perfect} For every $n \geq 3$, the set $\W_n$ does not contain any perfect lattices. Hence there are no extreme lattices in~$\W_n$ for~$n \geq 3$.
\end{cor}

There is another consequence of Theorem~\ref{SL_more} that we want to record: we also prove it in Section~\ref{SL}. Clearly, lattices in~$\W_n$ contain bases of minimal vectors. In fact, WR nearly orthogonal lattices satisfy a stronger property, which they have in common with such lattices as root lattices~$A_n$, for instance (see \textsection{4.2} of~\cite{martinet} for the definition and detailed properties of the $A_n$ family).

\begin{cor} \label{min_basis} Let $L \in \W^*_n$. Then any $n$ linearly independent vectors from~$S(L)$ form a basis for~$L$.
\end{cor}

Let us now look at the coherence of WR nearly orthogonal lattices in more details. Define a dimensional constant
\begin{equation}
\label{cn_def0}
c_n = \frac{\sqrt{(n-2)^2 + 16(n-1)} - (n-2)}{8(n-1)}.
\end{equation}
We prove the following result.

\begin{thm} \label{main_coh} The following statements are true.
\begin{enumerate}

\item For a lattice $L \in \W^*_n$, $C(L) = 1/2$ if and only if~$|S(L)| > 2n$.

\item Let $B = \{ \bb_1,\dots,\bb_n \}$ in $\real^n$ be a collection of linearly independent unit vectors such that
$$\max_{1 \leq i < j \leq n} |(\bb_i,\bb_j)| \leq c_n.$$
Then $L = \spn_{\zed} B$ is in~$\W^*_n$.

\end{enumerate}
\end{thm}

\noindent
We prove Theorem~\ref{main_coh} in Section~\ref{coherence}, where we also demonstrate a family of lattices~$A_n^*$ outside of~$\W_n$ with coherence~$1/n$, which tends to~$\W_n$ with respect to coherence as $n \to \infty$ in the sense that $\lim_{n \to \infty} (c_n / (1/n)) = 1$. This shows that~$c_n$ is asymptotically sharp, so WR lattices with coherence $< 1/n$ tend to be nearly orthogonal as $n \to \infty$. For comparison, Proposition~1.1 of~\cite{lat_ang} implies that the coherence of the set of the first $k$ shortest vectors of a random lattice in~$\real^n$ (where $k$ does not depend on $n$) tends to $O(1/\sqrt{n})$ as $n \to \infty$. In Section~\ref{coherence}, we also prove a result in the spirit of Diophantine approximation about an infinite family of integral WR lattices in the plane with coherence tending to~$0$. A summary of our contributions and a possible direction for future research are discussed Section~\ref{conclude}. We are now ready to proceed.
\bigskip

\section{Packing density}
\label{dense}

In this section we prove Theorem~\ref{dense_coh}. We start with several auxiliary lemmas.

\begin{lem} \label{short} Let $\bb_1, \bb_2 \in \real^n$ be nonzero vectors with the angle $\pi/3 \leq \theta_1 \leq 2 \pi/3$ between them. Then
\begin{equation}
\label{short1}
\min \{ \| \alpha \bb_1 + \beta \bb_2 \| : \alpha, \beta \in \zed \text{ not both } 0 \} = \min \{ \|\bb_1\|, \|\bb_2\| \}.
\end{equation}
If $\|\bb_1\| = \|\bb_2\|$ and $\alpha, \beta \neq 0$, then $\| \alpha \bb_1 + \beta \bb_2 \| = \|\bb_1\|$ if and only if either $\theta_1 = \pi/3$ and the pair $(\alpha, \beta) = \pm (1,-1)$, or $\theta_1 = 2\pi/3$ and the pair $(\alpha, \beta) = \pm (1,1)$.
\end{lem}

\proof
Assume $\min \{ \|\bb_1\|, \|\bb_2\| \} = \|\bb_1\|$. Let $\bx = \alpha \bb_1 + \beta \bb_2$ for some $\alpha, \beta \in \zed$, not both~$0$. If $\alpha = 0$ or $\beta = 0$, then we immediately have $\|\bx\| \geq \|\bb_1\|$. Suppose that $\alpha,\beta \neq 0$, then
\begin{eqnarray}
\label{short_eq1}
\|\bx\|^2 & = & \alpha^2 \|\bb_1\|^2 + 2 \alpha \beta \|\bb_1\| \|\bb_2\| \cos \theta_1 + \beta^2 \|\bb_2\|^2 \nonumber \\
& \geq & \alpha^2 \|\bb_1\|^2 - | \alpha \beta | \|\bb_1\| \|\bb_2\| + \beta^2 \|\bb_2\|^2 \nonumber \\
& = & ( |\alpha| \|\bb_1\| - |\beta| \|\bb_2\|)^2 + | \alpha \beta | \|\bb_1\| \|\bb_2\| \geq \|\bb_1\|^2
\end{eqnarray}
and
\begin{equation}
\label{short_eq2}
\|\bx\| = \|\bb_1\|\ \Longleftrightarrow\ \|\bb_2\| = \|\bb_1\|,\ \cos \theta_1 = \pm 1/2, \text{ and } \alpha \beta = \pm 1.
\end{equation}
This completes the proof.
\endproof

\begin{lem} \label{12} Let $\bb_1,\bb_2,\bb_3$ be nonzero non-coplanar vectors in~$\real^n$, $n \geq 3$. Let $\theta_{ij}$ be the angle between $\bb_i$ and $\bb_j$, $1 \leq i \neq j \leq 3$, and let $\xi$ be the angle between $\bz = \bb_1 + \bb_2$ and $\bb_3$. Suppose that
$$\cos \theta_{13} = \cos \theta_{23} = \alpha > 0.$$
Then $| \cos \xi | > \alpha$, and hence $\bb_3$ makes a smaller angle with the plane spanned by $\bb_1$ and $\bb_2$ than with each of them.
\end{lem}

\proof
Notice that $\cos \xi = \frac{(\bz, \bb_3)}{\|\bz\| \|\bb_3\|}$, where
$$(\bz, \bb_3) = (\bb_1,\bb_3) + (\bb_2,\bb_3) = \alpha \|\bb_3\| ( \|\bb_1\| + \|\bb_2\| ) > \alpha \|\bb_3\| \|\bb_1 + \bb_2\|.$$
Further, the angle $\bb_3$ makes with the plane spanned by $\bb_1$ and $\bb_2$ is $\leq \xi$, since the vector $\bz = \bb_1 + \bb_2$ lies in that plane. The conclusion follows.
\endproof

\begin{lem} \label{sub_WR} Let $L \in \W_n$ and 
$$B = \{ \bb_1,\dots,\bb_n \}$$
be a weakly nearly orthogonal basis for $L$. Then for each $1 \leq k \leq n$, the lattice
$$L_k = \spn_{\zed} \left\{ \bb_1,\dots,\bb_k \right\}$$
is also WR.
\end{lem}

\proof
Arguing toward a contradiction, suppose this is not true. Since $L$ is WR, there must exist some $1 < k < n$ such that $L_k$ is not WR, but $L_{k+1}$ is. This means that the number of linearly independent minimal vectors of $L_k$ is $1 \leq m < k$; call these vectors $\bx_1,\dots,\bx_m$. On the other hand, $L_k \subsetneq L_{k+1}$ and $L_{k+1}$ must contain~$k+1$ linearly independent minimal vectors. Now, every vector $\bwy \in L_{k+1}$ is of the form
$$\bwy = \alpha \bx + \beta \bb_{k+1}$$
for some $\bx \in L_k$ and $\alpha, \beta \in \zed$. Then, by Lemma~\ref{short},
$$|L_{k+1}| = \min \{ \|\bwy\| : \bwy \in L_{k+1} \setminus \{\bo\} \} = \min \{ |L_k|, \|\bb_{k+1}\| \}.$$
Suppose first that $\|\bb_{k+1}\| < |L_k|$. Then $\|\bwy\| = |L_{k+1}|$ if and only if $\bwy = \pm \bb_{k+1}$: since $k+1 \geq 2$, this contradicts $L_{k+1}$ being WR. Next assume that $|L_k| < \|\bb_{k+1}\|$. Then $\|\bwy\| = |L_{k+1}|$ if and only if $\bwy \in S(L_k)$. However, there are only $m < k$ such linearly independent vectors, again contradicting $L_{k+1}$ being WR. Hence the only remaining option is that
$$\|\bx_1\| = \dots = \|\bx_m\| = \|\bb_{k+1}\|.$$
In this case,
\begin{equation}
\label{slk}
S(L_k) \cup \{ \pm \bb_{k+1} \} \subseteq S(L_{k+1}),
\end{equation}
but only $m+1 < k+1$ of these vectors are linearly independent.  Additionally, some vectors of the form $\bx \pm \bb_{k+1}$ for some $\bx \in S(L_k)$ may be in $S(L_{k+1})$. However, they are linearly dependent with the vectors in~\eqref{slk}. This means that the number of linearly independent vectors in $S(L_{k+1})$ is $\leq m+1 < k+1$, contradicting the assumption that $L_{k+1}$ is WR. Hence every $L_k$ must be WR, and this completes the proof.
\endproof

\begin{lem} \label{near_WR} Let $L \in \W^*_n$ and $B$ a nearly orthogonal basis for $L$. Then $B \subseteq S(L)$.
\end{lem}

\proof
First suppose that $B$ is weakly $\theta$-orthogonal, where $\theta > \pi/3$. Then Corollary~1 of~\cite{baraniuk} guarantees that $S(L) = \pm B$. Assume then that for some $1 \leq i \leq n-1$, the angle $\theta_i$ between $\bb_i$ and subspace spanned by $\bb_1,\dots,\bb_{n-1}$ is equal to~$\pi/3$. 

We argue by induction on $n \geq 2$. If $n=2$, let $\bb_1, \bb_2$ be a nearly orthogonal basis for $L$; assume $\|\bb_1\| \leq \|\bb_2\|$. Then Theorem~1 of~\cite{baraniuk} guarantees that $|L| = \|\bb_1\|$. Let $\theta_1 = \pi/3$ be the angle between $\bb_1,\bb_2$. Since $L$ is WR, there must exist some
$$\bx = \alpha \bb_1 + \beta \bb_2 \in L$$
such that $\beta \neq 0$ and $\|\bx\| = \|\bb_1\|$. If $\alpha = 0$, then $\bx = \beta \bb_2$, and so we must have $\beta = \pm 1$ and $\|\bb_2\| = \|\bb_1\|$. If $\alpha \neq 0$, then $\|\bb_2\| = \|\bb_1\|$ by~\eqref{short_eq1} and~\eqref{short_eq2} above. In either case, $\bb_1,\bb_2 \in S(L)$.

Next assume the lemma is true in all dimensions $\leq n-1$. Let us prove it for dimension~$n$. Let 
$$B = \{ \bb_1,\dots,\bb_{n-1}, \bb_n \}$$
be a nearly orthogonal basis for a lattice $L \in \W^*_n$. Theorem~1 of~\cite{baraniuk} guarantees that at least one of these basis vectors is a shortest vector for $L$, and we can assume that $B$ is ordered so that it is~$\bb_1$. Let $L_{n-1} = \spn_{\zed} \{ \bb_1,\dots,\bb_{n-1} \}$, then Lemma~\ref{sub_WR} implies that $L_{n-1} \in \W^*_{n-1}$ and
\begin{equation}
\label{bb1}
|L_{n-1}| = \|\bb_1\| = |L|,
\end{equation}
and $\|\bb_n\| \geq |L|$. We should remark that when we write $\W^*_{n-1}$ here (and below), we are identifying $\real^{n-1}$ with $\spn_{\real} L_{n-1}$. By induction hypothesis, we have $\bb_1,\dots,\bb_{n-1} \in S(L_{n-1}) \subseteq S(L)$. Hence we only need to prove that $\bb_n \in S(L)$. Since $L$ is WR, there must exist some $\bwy \in S(L) \setminus L_{n-1}$. Again, $L = \spn_{\zed} \{ L_{n-1}, \bb_n \}$, so
$$\bwy = \alpha \bx + \beta \bb_n$$
for some $\bx \in L_{n-1}$ and $\alpha, \beta \in \zed$ with $\beta \neq 0$. By Lemma~\ref{short},
$$|L| = \|\bwy\|  = \min \{ |L_{n-1}|, \|\bb_n\| \}.$$
Combining this observation with~\eqref{bb1} we have, in particular
$$\min\{ \| \alpha \bb_1 + \beta \bb_n \| : \alpha, \beta \in \zed,\ \beta \neq 0 \} = \|\bb_1\|,$$
while $\|\bb_n\| \geq \|\bb_1\|$ and the angle between $\bb_1$ and $\bb_n$ is in the interval~$[\pi/3, 2\pi/3]$. Then~\eqref{short_eq1} and~\eqref{short_eq2} imply that~$\bb_n \in S(L)$, and we are done.
\endproof

We are now ready to prove the theorem.

\proof[Proof of Theorem~\ref{dense_coh}]
First notice that, by Lemma~\ref{near_WR}, $B \subseteq S(L)$. To prove parts (1) and (2) of the theorem, we argue by induction on $n \geq 2$. First suppose that $n=2$. Let $L \in \W^*_2$ and let $B = \{ \bb_1, \bb_2 \}$ be a nearly orthogonal basis for $L$. In this case $\mu(B) = \nu(B) = C(L)$, and assume that $0 < C(L) < 1/2$. Then
$$L = \spn_{\zed} \left\{ \bb_1,\bb_2 \right\}$$
with $|L| = \|\bb_1\| = \|\bb_2\| = 1$ and the angle $\theta_1$ between $\bb_1$ and $\bb_2$ lies in the interval $(\pi/3,\pi/2)$. The packing density of $L$ is
$$\delta(L) = \frac{\pi |L|^2}{\det L} = \frac{\pi}{\sin \theta_1}.$$
Let us write $U(\theta_1)$ for the counterclockwise rotation matrix by the angle~$\theta_1$:
$$U(\theta_1) = \begin{pmatrix} \cos \theta_1 & -\sin \theta_1 \\ \sin \theta_1 & \cos \theta_1 \end{pmatrix}.$$
Without loss of generality, we can assume that $\bb_2 = U(\theta_1) \bb_1$. For a given $\eps > 0$, let $\theta_1' \in [\pi/3, \theta_1)$ and $\theta_1'' \in (\theta_1,\pi/2]$ be such that
\begin{equation}
\begin{aligned}
\label{sin_cos}
& \left( \sin \theta_1 - \sin \theta_1' \right)^2 + \left( \cos \theta_1 - \cos \theta_1' \right)^2 < \eps, \\
& \left( \sin \theta_1 - \sin \theta_1'' \right)^2 + \left( \cos \theta_1 - \cos \theta_1'' \right)^2 < \eps.
\end{aligned}
\end{equation}
Then the lattices
$$L' = \spn_{\zed} \left\{ \bb_1, U(\theta_1') \bb_1 \right\},\ L'' = \spn_{\zed} \left\{ \bb_1, U(\theta_1'') \bb_1 \right\}$$
with nearly orthogonal bases $B' = \{ \bb_1, U(\theta_1') \bb_1 \}$, $B'' = \{ \bb_1, U(\theta_1'') \bb_1 \}$, respectively, are in $\W^*_2 \cap \BB_{\eps}(L)$, 
$$1/2 \geq \mu(B') > \mu(B) = \nu(B) > \nu(B'') \geq 0,$$
and
$$
\begin{aligned}
& \delta(L') = \frac{\pi}{\sin \theta'_1} = \frac{\pi}{\sqrt{1 - \mu(B')^2}} > \frac{\pi}{\sqrt{1 - \mu(B)^2}} = \frac{\pi}{\sin \theta_1} = \delta(L),\\
& \delta(L'') = \frac{\pi}{\sin \theta''_1} = \frac{\pi}{\sqrt{1 - \nu(B'')^2}} < \frac{\pi}{\sqrt{1 - \nu(B)^2}} = \frac{\pi}{\sin \theta_1} = \delta(L).
\end{aligned}
$$
Hence
$$\delta(L') = \frac
{\sqrt{1 - \mu(B)^2}}{\sqrt{1 - \mu(B')^2}}\ \delta(L),\ \delta(L'') = \frac{\sqrt{1 - \nu(B)^2}}{\sqrt{1 - \nu(B'')^2}}\ \delta(L).$$

Now suppose the statement is true in any dimension $\leq n-1$. Let us prove it for~$n$. We start with part (1), so let $L \in \W^*_n$ and let $B = \{ \bb_1,\dots,\bb_n \}$ be a nearly orthogonal basis for $L$ so that $\mu(B) < 1/2$. Then~$B \subseteq S(L)$ by Lemma~\ref{near_WR}. Let $L_{n-1} = \spn_{\zed} \{ \bb_1,\dots,\bb_{n-1} \}$, then by Lemmas~\ref{sub_WR} and~\ref{near_WR}, $L_{n-1} \in \W^*_{n-1}$ and $B_{n-1} = \{\bb_1,\dots,\bb_{n-1} \} \subseteq S(L_{n-1})$. Further, since we can reorder $B$ as we like, we can assume that $\mu(B_{n-1}) = \mu(B) < 1/2$, and thus we can apply induction hypothesis to~$L_{n-1}$. Then there exists $L'_{n-1} \in \W^*_{n-1} \cap \BB_{\eps}(L_{n-1})$ with the nearly orthogonal basis $B'_{n-1} = \{ \bb'_1,\dots,\bb'_{n-1} \}$ such that
$$\mu(B'_{n-1}) > \mu(B_{n-1}),\ \delta(L'_{n-1}) = \frac{\sqrt{1 - \mu(B_{n-1})^2}}{\sqrt{1 - \mu(B'_{n-1})^2}}\ \delta(L_{n-1}).$$
Since we agreed to pick representatives of similarity classes that have minimal norm~$1$, we have
$$|L'_{n-1}| = \|\bb'_1\| = \dots = \|\bb'_{n-1}\| = |L_{n-1}| = |\bb_n\| = |L| = 1,$$
by Lemma~\ref{near_WR}. Now, let $L' = \spn_{\zed} \{ L'_{n-1}, \bb_n \}$. Since $L'_{n-1} \subset \spn_{\real} L_{n-1}$, $\bb_n$ makes the same angle $\theta_{n-1}$ with $\spn_{\real} \{ \bb'_1,\dots,\bb'_{n-1} \}$ and with $\spn_{\real} \{ \bb_1,\dots,\bb_{n-1} \}$, and so $L' \in \W^*_n$ with near orthogonal basis $B' = \{ \bb'_1,\dots,\bb'_{n-1},\bb_n \}$ and $|L'| = |L|$,
$$\mu(B') = \mu(B'_{n-1}) > \mu(B_{n-1}) = \mu(B),$$
and so $\delta(L') = \frac{\sqrt{1 - \mu(B)^2}}{\sqrt{1 - \mu(B')^2}}\ \delta(L)$, since
$$\det L' = (\det L'_{n-1}) \|\bb_n\| \sin \theta_{n-1},\ \det L'= (\det L_{n-1}) \|\bb_n\| \sin \theta_{n-1}.$$
This completes the proof of (1). The proof of (2) is completely analogous with $\mu$ replaced by $\nu$ and the corresponding inequalities reversed.

To prove part (3), assume $n \geq 3$ and suppose $\mu(B) = 1/2$ for some $L \in \W_n^*$ with a nearly orthogonal basis~$B = \{ \bb_1,\dots,\bb_n \}$. This means that all the angles between these basis vectors are equal to~$\pi/3$ or $2\pi/3$. In particular, $\bb_3$ makes such an angle with $\bb_1$ and $\bb_2$. But then Lemma~\ref{12} implies that $\bb_3$ makes an angle $< \pi/3$ with the plane spanned by $\bb_1,\bb_2$, contradicting near-orthogonality of the basis~$B$. Therefore we must have $\mu(B) < 1/2$ when $n \geq 3$. On the other hand, if $n=2$ and $L \in \W_2^*$, then for any nearly orthogonal basis $B$, $\mu(B) \leq C(L)$ is simply the cosine of the angle between the minimal basis vectors, which is in the interval $[\pi/3,\pi/2]$ and is equal to $\pi/3$ precisely in the case of the hexagonal lattice $\begin{pmatrix} 1 & 1/2 \\ 0 & \sqrt{3}/2 \end{pmatrix} \zed^2$.

Finally, for part (4) notice that $C(L) = 0$ if and only if all the angles between minimal vectors are equal to~$\pi/2$, which happens precisely in the case of the integer lattice~$\zed^n$.
\endproof
\bigskip

\section{Minimal vectors}
\label{SL}

In this section we construct families of latices in $\W_n$ with many minimal vectors, but also prove that for any~$L \in \W_n$, $|S(L)| \leq 4n-2$. This will establish Theorem~\ref{SL_more}. We begin with two constructions. We write $A_2$ for the $2$-dimensional root lattice (isometric to the hexagonal lattice), normalized to have minimal norm~$1$. We also write $\oplus$ for the orthogonal direct sum of lattices.

\begin{lem} \label{C1} Let $n \geq 2$. For each $0 \leq m \leq n/2$, let
\begin{equation}
\label{lmn}
L_{n,m} = \bigoplus_{i=1}^m A_2 \bigoplus_{j=1}^{n-2m} \zed,
\end{equation}
where a direct sum is taken to be empty if the upper limit on the index is~$0$. This is a strongly eutactic lattice contained in~$\W^*_n$ with
$$|S(L_{n,m})| = 2m + 2n.$$
If we take $m = [n/2]$, we have
$$|S(L_{n,m})| = \left\{ \begin{array}{ll}
3n & \mbox{if $n$ is even} \\
3n-1 & \mbox{if $n$ is odd}.
\end{array}
\right.$$
\end{lem}

\proof
Let $B = \{ \bb_1,\dots,\bb_n \} \subset \real^n$ be a basis of unit vectors satisfying the following condition:
\begin{eqnarray}
\label{condition}
& \text{For all } 1 \leq i \leq n \text{ there exists at most one } 1 \leq t \leq n \text{ such that:} \nonumber \\
&  (\bb_i,\bb_t) = 1/2 \text{ and } \text{for all } 1 \leq k \leq n, k \neq i,t,\ (\bb_i,\bb_k) = 0.
\end{eqnarray}
Notice that such a basis is nearly orthogonal. Indeed, let $1 \leq i \leq n$ and let $V$ be a subspace of~$\real^n$ spanned by some of the other vectors of~$B$, say 
$$V = \spn_{\real} \{ \bb_{j_1},\dots,\bb_{j_m} : 1 \leq j_1,\dots,j_m \leq n, j_k \neq i\ \text{for all } k = 1,\dots,m \}.$$
Let 
$$\bx = \sum_{k=1}^m a_k \bb_{j_k} \in V$$
be a unit vector for some real coefficients $-1 \leq a_1,\dots,a_m \leq 1$. Since $\bb_i$ is orthogonal to every other vector of $B$ except for (possibly) $\bb_t$ and $(\bb_i,\bb_t) = 1/2$, we have $|(\bx,\bb_i)| \leq 1/2$. Thus $\bb_i$ makes an angle $\geq \pi/3$ with each such subspace~$V$. Thus if $L = \spn_{\zed} B$, then $L \in \W_n^*$.

Now, let $0 \leq m \leq n/2$ and let~$L_{n,m}$ be as in~\eqref{lmn}. Then~$L_{n,m}$ is spanned over~$\zed$ by a unit basis
$$\bb_{11},\bb_{12},\dots,\bb_{m1},\bb_{m2},\bc_1,\dots,\bc_{n-2m},$$
where each pair $\bb_{i1},\bb_{i2}$ spans a copy of~$A_2$ and hence has inner product~$1/2$, $\bc_j$'s span~$\zed^{n-2m}$, and hence are orthogonal to each other, and each pair $\bb_{i1},\bb_{i2}$ is orthogonal to every other pair and to all~$\bc_j$'s. Hence this basis satisfies condition~\eqref{condition}, therefore~$L_{n,m} \in \W^*_n$. 

Let us count the number of minimal vectors in~$L_{n,m}$. Notice that each of the~$m$ copies of~$A_2$ in the orthogonal direct sum contributes~$6$ minimal vectors:
$$\pm \bb_{i1}, \pm \bb_{i2}, \pm (\bb_{i1} - \bb_{i2}),$$
and each of the~$n-2m$ copies of~$\zed$ contributes two minimal vectors: $\pm \bc_j$. There are no other minimal vectors. Hence we have
$$|S(L_{m,n})| = 6m + 2(n-2m) = 2m + 2n.$$

Finally notice that $L_{n,m}$ is an orthogonal direct sum of strongly eutactic lattices with equal minimal norm. Therefore it is strongly eutactic by Theorem~3.6.13 of~\cite{martinet}. This completes the proof of the lemma.
\endproof

\begin{lem} \label{C2} For every $n \geq 2$ there exist $L \in \W_n$ such that $|S(L)| \geq 4n-2$. 
\end{lem}

\proof
To prove this lemma, we demonstrate another construction. Let $\bb_1,\bb_2 \in \real^n$ be unit vectors with angle $\pi/3$ between them, and let $\Pi_1$ be the plane spanned by them. Then $\bb_1 - \bb_2$ is also a unit vector in~$\Pi_1$ by Lemma~\ref{short}. Let $\bb_3$ be a unit vector in~$\real^n$ making an angle $\pi/3$ with $\Pi_1$ and with $\bb_1-\bb_2$, i.e., the orthogonal projection of~$\bb_3$ onto~$\Pi_1$ is along the line spanned by~$\bb_1-\bb_2$. Then $\bb_1 - \bb_2 - \bb_3$ is also a unit vector. Let $\bb_4$ be a unit vector in~$\real^n$ making an angle of $\pi/3$ with~$\bb_1-\bb_2-\bb_3$ and with the $3$-dimensional subspace spanned by $\bb_1,\bb_2,\bb_3$. Once again, $\bb_1-\bb_2-\bb_3-\bb_4$ is also a unit vector. Continuing in the same manner, we construct a basis $B = \{ \bb_1,\dots,\bb_n \}$ and take $L = \spn_{\zed} B$. Then $L \in \W_n$ by construction. Further, for $n \geq 2$ the vectors 
$$\pm \bb_k \text{ for all } 1 \leq k \leq n,\ \pm \left( \bb_1 - \sum_{i=2}^k \bb_i \right)\ \text{for all } 2 \leq k \leq n$$
are contained in~$S(L)$. The number of these vectors is~$4n-2$.
\endproof

\begin{ex} \label{ex_34} Let us show examples of our constructions in the proofs of Lemmas~\ref{C1} and~\ref{C2} above for $n=3,4$. For an example of the first construction when~$n=3$, we can take
$$L_1 = \spn_{\zed} \left\{ \frac{1}{\sqrt{2}} \begin{pmatrix} 1 \\ 1 \\ 0 \end{pmatrix}, \frac{1}{\sqrt{2}} \begin{pmatrix} 1 \\ 0 \\ 1 \end{pmatrix}, \frac{1}{\sqrt{3}} \begin{pmatrix} -1 \\ 1 \\ 1 \end{pmatrix} \right\},$$
which is a lattice in~$\W_3^*$ with~$8$ minimal vectors. For $n=4$, take
$$L_2 = \spn_{\zed} \left\{ \frac{1}{\sqrt{2}} \begin{pmatrix} 1 \\ 1 \\ 0 \\ 0 \end{pmatrix}, \frac{1}{\sqrt{2}} \begin{pmatrix} 1 \\ 0 \\ 1 \\ 0 \end{pmatrix}, \frac{1}{2} \begin{pmatrix} -1 \\ 1 \\ 1 \\ 1 \end{pmatrix}, \frac{1}{2} \begin{pmatrix} 1 \\ -1 \\ -1 \\ 1 \end{pmatrix} \right\},$$
which is a lattice in~$\W_4^*$ with~$12$ minimal vectors. The presented bases for these lattices satisfy~\eqref{condition}.

Here is also a~$3$-dimensional example of the second construction:
$$L_3 = \spn_{\zed} \left\{ \frac{1}{\sqrt{2}} \begin{pmatrix} 1 \\ 1 \\ 0 \end{pmatrix}, \frac{1}{\sqrt{2}} \begin{pmatrix} 1 \\ 0 \\ 1 \end{pmatrix}, \frac{1}{2 \sqrt{2}} \begin{pmatrix} \sqrt{2} \\ 1-\sqrt{2} \\ -(1+\sqrt{2}) \end{pmatrix} \right\}.$$
This is a lattice in~$\W_3$ with~$10$ minimal vectors, however it is not in~$\W^*_3$, since the ordering of the minimal basis
$$\left\{ \frac{1}{\sqrt{2}} \begin{pmatrix} 1 \\ 1 \\ 0 \end{pmatrix}, \frac{1}{2 \sqrt{2}} \begin{pmatrix} \sqrt{2} \\ 1-\sqrt{2} \\ -(1+\sqrt{2}) \end{pmatrix}, \frac{1}{\sqrt{2}} \begin{pmatrix} 1 \\ 0 \\ 1 \end{pmatrix} \right\}$$
is not weakly nearly orthogonal: indeed, cosine of the angle between the plane spanned by the first two of these vectors and the third one is~$\sqrt{2/5} > 1/2$. In fact, since every $3$-dimensional lattice constructed as in the proof of Lemma~\ref{C2} is isometric to~$L_3$, all of them would be in~$\W_3$, but not in~$\W_3^*$. Furthermore, this implies that construction of Lemma~\ref{C2} never produces lattices in~$\W_n^*$ for $n \geq 3$: just reorder the first three vectors as in the example~$L_3$.
\end{ex}

\begin{lem} \label{4n-2} For any $L \in \W_n$, $|S(L)| \leq 4n-2$. 
\end{lem}

\proof
We argue by induction on~$n$. If $n=2$, the hexagonal lattice has largest set of minimal vectors, which has cardinality $6 = 4 \times 2 -2$. Now assume the statement is true in all dimensions~$\leq n-1$. We prove it in dimension~$n > 2$. Let $L \in \W_n$, and let $B = \{ \bb_1,\dots,\bb_n\}$ be a weakly nearly orthogonal basis for~$L$. Then the lattice $L' = \spn_{\zed} \{ \bb_1,\dots,\bb_{n-1} \} \in \W_{n-1}$ by Lemma~\ref{sub_WR}, and hence has at most $4(n-1) - 2 = 4n-6$ minimal vectors by induction hypothesis. Suppose that $\bwy \in S(L)$ is not contained in~$L'$. Then either $\bwy = \bb_n$ or
$$\bwy = \alpha \bx + \beta \bb_n$$
for some $\bx \in L'$, $\beta \neq 0$, and
$$\|\bwy\| = |L'| = |L| = \|\bb_n\|.$$
If $\bwy \neq \bb_n$, Lemma~\ref{short} implies that $\|\bx\| = \|\bb_n\|$ and the angle between~$\bx$ and~$\bb_n$ is~$\pi/3$. By Lemma~\ref{12}, there can exist no more than one vector in~$L'$ with which~$\bb_n$ makes an angle~$\pi/3$: otherwise it would make an angle $< \pi/3$ with the subspace $\spn_{\real} L'$. Hence the total number of minimal vectors of~$L$ which are outside of~$L'$ is no greater than~$4$, and so
$$|S(L)| \leq |S(L')| + 4 \leq 4n-2.$$
\endproof

\begin{cor} \label{C2-1} For every $n \geq 3$ and $m \geq 1$ such that
\begin{equation}
\label{c211}
m \leq \left\{ \begin{array}{ll}
\frac{n-2}{2} & \mbox{if $n$ is even} \\
\frac{n-1}{2} & \mbox{if $n$ is odd},
\end{array}
\right.
\end{equation}
there exist $L \in \W_n$ such that 
\begin{equation}
\label{c212}
|S(L)| = \left\{ \begin{array}{ll}
3n+2m & \mbox{if $n$ is even} \\
3n-1+2m & \mbox{if $n$ is odd}.
\end{array}
\right.
\end{equation}
\end{cor}

\proof
We argue by induction on~$n \geq 3$. If $n=3$, we must have $m=1$, and so 
$$3n-1+2m = 9 - 1 + 2 = 10 = 4n -2.$$
If $n=4$, again we have $m=1$, so
$$3n+2m = 12 + 2 = 14 = 4n -2.$$
Hence the existence of such a lattice $L$ in $\W_3$ or $\W_4$ follows directly from Lemma~\ref{C2}. Assume then that $n > 4$, and the result holds in all dimensions~$\leq n-1$. Let us prove it for~$n$. First notice that if there is equality in~\eqref{c211}, then the result again follows from Lemma~\ref{C2}. Hence let us assume that $m$ is strictly less than the right hand side of~\eqref{c211}, that is
$$m \leq \left\{ \begin{array}{ll}
\frac{n-4}{2} = \frac{(n-2) - 2}{2} & \mbox{if $n-2$ is even} \\
\frac{n-3}{2} = \frac{(n-2) - 1}{2} & \mbox{if $n-2$ is odd}.
\end{array}
\right.$$
By the induction hypothesis, there must exist~$L' \in \W_{n-2}$ with
$$|S(L')| = \left\{ \begin{array}{ll}
3(n-2)+2m & \mbox{if $n-2$ is even} \\
3(n-2)-1+2m & \mbox{if $n-2$ is odd}.
\end{array}
\right.$$
Let us embed~$L'$ into the $(n-2)$-dimensional coordinate subspace in~$\real^n$ corresponding to the last two coordinates being~$0$ and let~$W$ be the orthogonal complement of this subspace in~$\real^n$. Let $L''$ be a lattice of unit minimal norm in~$W$, similar to $A_2$, i.e., spanned by a pair of unit vectors~$\bc_1,\bc_2$ with angle~$\pi/3$ between them. Then~$S(L'') = \{ \pm \bc_1, \pm \bc_2, \pm (\bc_1 - \bc_2) \}$. Now, let $L = L' \oplus L''$. Then $S(L) = S(L') \cup S(L'')$, and so $|S(L)| = |S(L')| + 6$. This establishes~\eqref{c212}.
\endproof

\begin{ex} \label{explicit_5d} Here we present an explicit construction of the first case where there exists a lattice between the two extreme ends of Lemmas~\ref{C1} and~\ref{C2} (when $n=5$). First, we construct a 4-dimensional lattice $L_4 \in \W_4$ with $|S(L)|=14$ by following the procedure in Lemma~\ref{C2}. This gives
$$L_4 = \spn_{\zed} \left\{ \begin{pmatrix} 1 \\ 0 \\ 0 \\ 0 \end{pmatrix}, \frac{1}{2} \begin{pmatrix} 1 \\ \sqrt{3} \\ 0 \\ 0 \end{pmatrix}, \frac{1}{4} \begin{pmatrix} 1 \\ -\sqrt{3} \\ 2\sqrt{3} \\ 0 \end{pmatrix},\frac{1}{8} \begin{pmatrix} 1 \\ -\sqrt{3} \\ -2\sqrt{3} \\ 4\sqrt{3} \end{pmatrix}  \right\}.$$
Then, to produce a 5-dimensional lattice $L_5 \in \W_5$ with $|S(L)|=16=3n+1$, we need simply to add a copy of $\zed$ orthogonal to the subspace of $\real^5$ spanned by the basis of $L_4$. This gives us
$$L_5 = \spn_{\zed} \left\{ \begin{pmatrix} 1 \\ 0 \\ 0 \\ 0 \\ 0 \end{pmatrix}, \frac{1}{2} \begin{pmatrix} 1 \\ \sqrt{3} \\ 0 \\ 0 \\ 0 \end{pmatrix}, \frac{1}{4} \begin{pmatrix} 1 \\ -\sqrt{3} \\ 2\sqrt{3} \\ 0 \\ 0 \end{pmatrix},\frac{1}{8} \begin{pmatrix} 1 \\ -\sqrt{3} \\ -2\sqrt{3} \\ 4\sqrt{3} \\ 0 \end{pmatrix}, \begin{pmatrix} 0 \\ 0 \\ 0 \\ 0 \\ 1 \end{pmatrix}  \right\}.$$
By the same argument as at the end of Example~\ref{ex_34}, lattices produced using Corollary~\ref{C2-1} are again in~$\W_n$, but not in~$\W_n^*$, since they still follow the construction of Lemma~\ref{C2}.
\end{ex}
\smallskip

The strongly eutactic lattices constructed in Lemma~\ref{C1} are orthogonal direct sums of copies of~$A_2$ and~$\zed$, which is not so surprising. It is more interesting that~$\W_n$ can contain irreducible eutactic lattices which are also not in~$\W^*_n$; we now demonstrate such an example for~$n=3$.

\begin{ex} \label{eutactic_3d} Let us consider eutactic lattices in dimensions~$2$ and~$3$. In~$\real^2$, there are only two eutactic lattices: $\zed^2$ and $A_2$, and both of them are in~$\W^*_2$ by Lemma~\ref{C1}. Let us write $\perp$  for the orthogonal direct sum of lattices, then in dimension~$3$ there are five eutactic lattices: $\zed^3$, $A_2 \perp \zed$, $A_3$, $A_3^*$ and $K_3'$ (see Example 9.5.1 (6) on p.~345 of~\cite{martinet}). The lattices $\zed^3$ and $A_2 \perp \zed$ are in~$\W^*_3$ by Lemma~\ref{C1}. The root lattice $A_3$ has~$12$ minimal vectors, and hence is not in~$\W_3$ by Theorem~\ref{SL_more}. The lattice~$A_3^*$ is also not in~$\W_3$ (see Example~\ref{ETF} below). The one remaining lattice is $K'_3$ (see Section~8.5 of~\cite{martinet} for its construction), which is spanned by a unit basis
$$\bb_1 = \begin{pmatrix} 1 \\ 0 \\ 0 \end{pmatrix},\ \bb_2 = \begin{pmatrix} -1/2 \\ \sqrt{3}/2 \\ 0 \end{pmatrix},\ \bb_3 = \begin{pmatrix} -1/2 \\ 0 \\ \sqrt{3}/2 \end{pmatrix}.$$
This lattice has~$10$ minimal vectors: $\pm \bb_i$ for $1 \leq i \leq 3$, $\pm (\bb_1 + \bb_2)$, $\pm (\bb_1 + \bb_3)$. It is irreducible, eutactic, but not perfect, and not strongly eutactic. We will show that~$K'_3 \in \W_3$, but~$K'_3 \not\in \W^*_3$. Indeed, let $\theta_{ij}$ be the angle between~$\bb_i$ and $\bb_j$, $1 \leq i < j \leq 3$. Also define
$$\Pi_1 = \spn_{\real} \{ \bb_1,\bb_2 \},\ \Pi_2 = \spn_{\real} \{ \bb_2, \bb_3 \},$$
and let $\nu_1$ be the angle between $\bb_3$ and $\Pi_1$, $\nu_2$ the angle between $\bb_1$ and $\Pi_2$. It is then easy to check that
$$| \cos \theta_{12} | = | \cos \theta_{13} | = \frac{1}{2},\ | \cos \theta_{23} | = \frac{1}{4}.$$
Also, $| \cos \nu_1 | = \frac{1}{2}$, and hence $B = \{ \bb_1, \bb_2, \bb_3 \}$ is a weakly nearly orthogonal basis, so~$K'_3 \in \W_3$. On the other hand, let 
$$\bx = \frac{1}{2} ( \bb_2 + \bb_3 ) = \begin{pmatrix} -1/2 \\ \sqrt{3}/4 \\ \sqrt{3}/4 \end{pmatrix} \in \Pi_2,$$ 
and let $\mu$ be the angle between~$\bb_1$ and~$\bx$. Notice that
$$| \cos \nu_2 | \geq | \cos \mu | = \sqrt{ \frac{2}{5}} > \frac{1}{2}.$$
This means that the ordering of the basis $\{ \bb_2, \bb_3, \bb_1 \}$ is not weakly nearly orthogonal, and hence~$K'_3 \not\in \W^*_3$. Notice also that~$K'_3$ is not similar to the lattice~$L_3$ in Example~\ref{ex_34}.
\end{ex}

\begin{rem} \label{eut_min} A theorem of A. Ash~\cite{ash1} (see also~\cite{ash2}) asserts that all the critical points of the packing density function~$\delta$ occur at eutactic lattices. By Voronoi's theorem, we know that these are maxima if and only if the corresponding eutactic lattice is also perfect. Non-perfect eutactic lattices may or may not be minima: combining our observations on eutactic lattices in~$\W_n$ with our Theorem~\ref{dense_coh}, we see that many of them are not minima, not even among well-rounded lattices. On the other hand, two lattices $L_1$ and $L_2$ are said to be in the same {\it minimal class} if there exists~$U \in \GL_n(\real)$ such that~$L_2 = UL_1$ and $S(L_2) = US(L_1)$. Theorem~9.4.1 of~\cite{martinet} asserts that~$\delta$ attains its minimum on a given minimal class at a weakly eutactic lattice, if one exists (each minimal class has no more than one weakly eutactic lattice). 
\end{rem}

Finally we show that lattices in~$\W^*_n$ cannot have too many minimal vectors, consistent with examples~$L_3$, $L_4$, $L_5$, and $K'_3$ demonstrated above.

\begin{lem} \label{notstar} Let $n \geq 3$ and $L \in \W_n$ be such that $|S(L)| > 3n$. Then $L \not\in \W^*_n$.
\end{lem}

\proof
Let $L \in \W^*_n$ and let $B = \{ \bb_1, \dots, \bb_n \}$ be its nearly orthogonal basis. We want to prove that $|S(L)| \leq 3n$. Suppose that for some~$1 \leq k \leq n$
$$\bx_1 = \alpha_1 \bb_k + \bwy \in S(L),\ \bx_2 = \alpha_2 \bb_k + \bz \in S(L),$$
where $0 \neq \alpha_1,\alpha_2 \in \zed$ and $\bo \neq \bwy, \bz \in \spn_{\zed} B \setminus \{ \bb_k \}$. Then by Lemma~\ref{short}, we must have $\alpha_1,\alpha_2 = \pm 1$, $\bwy, \bz \in S(L)$, and the angles between $\bb_k$ and $\bwy$, $\bz$ equal to $\pi/3$ or $2\pi/3$. In this case Lemma~\ref{12} implies that the angle~$\bb_k$ makes with the space spanned by the rest of vectors of $B$ is less than~$\pi/3$, which contradicts the assumption that $L$ is in~$\W^*_n$. Hence there can be at most one $\pm$ pair of vectors in $S(L)$ besides $\pm \bb_k$ which is expressible as an integral linear combination of the vectors of $B$ with a nonzero coefficient in front of~$\bb_k$, and this is true for every~$1 \leq k \leq n$. Thus a maximal possible number of minimal vectors for $L$ is achieved by the construction described in Lemma~\ref{C1} with $3n$ or $3n-1$ vectors, depending on whether~$n$ is even or odd. 
\endproof
\smallskip

\proof[Proof of Theorem~\ref{SL_more}]
The theorem now follows upon combining Lemmas~\ref{C1}, \ref{C2}, \ref{4n-2} and~\ref{notstar} with Corollary~\ref{C2-1}.
\endproof

\proof[Proof of Corollary~\ref{min_basis}]
We argue by induction on~$n \geq 2$. If $n=2$, then $|S(L)| = 4$ unless $L$ is the hexagonal lattice, in which case it has $6$ minimal vectors. In either case, the result is immediate by direct verification. Suppose now the result is established in all dimensions $\leq n-1$. Let us prove it for~$n$. By Lemma~\ref{notstar}, $|S(L)| \leq 3n$, out of which $n$ pairs of vectors are the nearly orthogonal basis vectors $\pm B = \pm \{ \bb_1,\dots,\bb_n \}$. Let $X = \{ \bx_1,\dots,\bx_n\} \subset S(L)$ be any $n$ linearly independent vectors. Then at least $n - [n/2]$ of them are vectors from~$\pm B$. Let $\bb_k$ be one of these vectors. There is at most one other vector, say $\bx_1 \in X$, which is a linear combination of some~$\bb_i$'s with $\pm 1$ coefficients and a nonzero coefficient in front of~$\bb_k$: we can write this $\bx_1$ as $\bx'_1 \pm \bb_k$. Then
$$\spn_{\zed} \{ \bx_1, \dots, \bx_n \} = \spn_{\zed} \{ \bx'_1, \dots, \bx_n \},$$
and $\spn_{\zed} (X \setminus \{ \bb_k \}) \subseteq L'_k := \spn_{\zed} ( B \setminus \{ \bb_k \})$ with $X \setminus \{ \bb_k \} \subset S(L'_k)$. Applying the induction hypothesis to $L'_k$, we see that $X \setminus \{ \bb_k \}$ is a basis for $L'_k$. Since $L = \spn_{\zed} \{ L'_k, \bb_k \}$, we conclude that~$X$ is a basis for~$L$.
\endproof

\bigskip

\section{Coherence}
\label{coherence}

We now discuss the coherence of lattices in some more details. As indicated in~\cite{ejc}, one might expect that many extreme lattices have coherence~$=1/2$. Certainly this is true for the standard root lattices~$A_n$, $D_n$, $E_6$, $E_7$ and $E_8$, as witnessed by the Coxeter-Dynkin diagrams (see, for instance, Theorem~4.6.3 of~\cite{martinet}). On the other hand, there are also extreme lattices that have coherence less than~$1/2$. Consider, for instance, the Coxeter-Barnes lattice~$A_n^r$, which is best defined as a lattice of rank~$n$ in~$\real^{n+1}$ spanned over~$\zed$ by the basis
$$\be_1 - \be_2,\dots, \be_1 - \be_n,\ \frac{1}{r} \left( n\be_1 - \sum_{i=2}^n \be_i \right),$$
where $\be_i$ are standard basis vectors in~$\real^{n+1}$, $n \geq 7$ and $1 < r < n+1$ is a divisor of~$n+1$. With parameters as specified, these lattices are known to be perfect and strongly eutactic, hence extreme (see Theorem~5.2.1 of~\cite{martinet}). If~$r=(n+1)/2$, these lattices have coherence~$< 1/2$ (see Proposition~5.2.3 of~\cite{martinet}).

These considerations raise a question: what is the coherence of a WR nearly orthogonal lattice? Well, it can be~$1/2$, as in the constructions Lemmas~\ref{C1} and~\ref{C2} above. In fact, this is the case for any $L \in \W^*_n$ with $|S(L)| > 2n$. The following proposition is part (1) of Theorem~\ref{main_coh}.

\begin{prop} \label{SL_coh} Let $L \in \W^*_n$. Then $C(L) = 1/2$ if and only if~$|S(L)| > 2n$.
\end{prop}

\proof
Let $L \in \W^*_n$ and $B$ be a weakly nearly orthogonal basis for~$L$. Suppose $|S(L)| = 2n$. Then $S(L) = \pm B$ by Lemma~\ref{near_WR}. Assume $C(L) = 1/2$. Then there are some two vectors, say, $\bb_i, \bb_j \in S(L)$ so that the angle between them is~$\pi/3$ or $2\pi/3$. Lemma~\ref{short} then implies that one of the vectors~$\bb_i \pm \bb_j$ is also in~$S(L)$, contradicting the fact that $S(L) = \pm B$. Hence, if $|S(L)| = 2n$, coherence must be $< 1/2$.

Now suppose $|S(L)| > 2n$. Then $\pm B \subsetneq S(L)$, so there must exist some $\bx \in S(L) \setminus \pm B$. Suppose that this $\bx$ is a linear combination of some $m \geq 2$ vectors of~$B$, say
$$\bx = \sum_{k=1}^m \alpha_k \bb_{i_k},$$
where $2 \leq m \leq n$, $1 \leq i_1 < \dots < i_m \leq n$, $\alpha_1,\dots,\alpha_m \in \zed$. We will prove that $C(L) = 1/2$. If $m=2$, then $\bx = \alpha_1 \bb_{i_1} + \alpha_2 \bb_{i_2} \in S(L)$ and
$$\|\bx\| = \|\bb_{i_1}\| = \|\bb_{i_2}\|.$$
Lemma~\ref{short} then implies that $\alpha_1, \alpha_2 = \pm 1$ and the angle between $\bb_{i_1}$ and $\bb_{i_2}$ is $\pi/3$ or $2\pi/3$. This implies that $C(L) = 1/2$. Assume now $m > 2$. Let $\bwy = \sum_{k=1}^{m-1} \alpha_k \bb_{i_k}$, then $\bx = \bwy + \alpha_m \bb_{i_m}$. Since $\bwy \in \spn_{\real} \{ \bb_{i_1},\dots,\bb_{i_{m-1}} \}$, the angle $\theta$ between $\bwy$ and $\bb_{i_m}$ is in the interval $[\pi/3,2\pi/3]$. Then Lemma~\ref{short} implies that
$$\|\bx\| \geq \min \{ \|\bwy\|, \|\bb_{i_m}\| \} = \|\bb_{i_m}\| = \|\bx\|,$$
which is only possible if $\bwy \in S(L)$, $\alpha_m = \pm 1$ and $\theta = \frac{\pi}{3}$ or $\frac{2\pi}{3}$. Hence $C(L)= \frac{1}{2}$.
\endproof

On the other hand, $\zed^n \in \W^*_n$ and it is reasonable to expect that lattices with very low coherence should be in~$\W_n^*$: large angles between minimal vectors make large angles in the sequence defining near orthogonality more likely. One can then ask how low is low enough? In other words, does there exist some dimensional constant~$c_n$ such that whenever~$C(L) \leq c_n$, the lattice~$L$ is necessarily in~$\W^*_n$?

\begin{ex} \label{ETF} Define a cyclic frame
$$\bb_1 := \frac{1}{\sqrt{n^2+n}} \left(\begin{array}{r} -n \\ 1 \\ \vdots \\ 1  \end{array}\right), \ldots, \bb_n := \frac{1}{\sqrt{n^2+n}} \left(\begin{array}{r} 1\\ \vdots \\-n \\ 1 \end{array}\right),$$
and
$$\bb_{n+1} := \frac{1}{\sqrt{n^2+n}} \left(\begin{array}{r} 1\\ 1\\ \vdots \\-n \end{array}\right) = -(\bb_1 + \dots + \bb_n).$$
Let $L = \spn_{\zed} \{ \bb_1,\dots,\bb_n \}$, then $L$ is a full-rank lattice in the hyperplane
$$H = \left\{ \bx \in \real^{n+1} : \sum_{i=1}^{n+1} x_i = 0 \right\},$$
and $S(L) = \{ \bb_1,\dots,\bb_n,\bb_{n+1}\}$. Identifying $\real^n$ with $H$, we view the root lattice $A_n$ as $\zed^{n+1} \cap H$, then $L$ is similar to the lattice~$A_n^*$, the dual of $A_n$, i.e.,
$$A^*_n := \left\{ \bx \in H : (\bx, \bwy) \in \zed\ \text{for all } \bwy \in A_n \right\}.$$
Then $C(L) = 1/n$, $|S(L)| = 2n + 2$; see~\cite{ejc}, \cite{etf}, \cite{martinet} for further details on this lattice. On the other hand, $L$ is not contained in~$\W_n$ by Proposition~\ref{SL_coh}, since $C(L) < 1/2$ while $|S(L)| > 2n$.
\end{ex}

\noindent
Thus Example~\ref{ETF} shows that a WR lattice with coherence even as low as~$1/n$ still does not have to be nearly orthogonal, i.e., $c_n < 1/n$. This being said, we can prove the following criterion, which is part~(2) of Theorem~\ref{main_coh}.

\begin{prop} \label{coh_bnd} Let $B = \{ \bb_1,\dots,\bb_n \}$ in $\real^n$ be a collection of linearly independent unit vectors such that
\begin{equation}
\label{cn_def}
\max_{1 \leq i < j \leq n} |(\bb_i,\bb_j)| \leq c_n,
\end{equation}
where $c_n$ is as in~\eqref{cn_def0}. Then $L = \spn_{\zed} B$ is in~$\W^*_n$.
\end{prop}

\proof
Let us prove that if~\eqref{cn_def} holds, then $B$ is a nearly orthogonal basis for $L = \spn_{\zed} B$. Without loss of generality, assume that $B = \{ \bb_1,\dots,\bb_n \}$ is an arbitrary ordering of~$B$. Then we only need to prove that~\eqref{cn_def} forces the angle $\theta$ between $\Pi_{k-1} := \spn_{\real} \{ \bb_1,\dots,\bb_{k-1} \}$ and $\bb_k$ to be $\geq \pi/3$ for every $2 \leq k \leq n$. Let $\bx \in \Pi_{k-1}$ be a vector so that
$$\aa(\bx,\bb_k) = \theta.$$
Let us write $\bx = \sum_{i=1}^{k-1} \alpha_i \bb_i$ for some $\alpha_1,\dots,\alpha_{k-1} \in \real$. Then
$$\|\bx\|^2 = \sum_{i,j=1}^{k-1} \alpha_i \alpha_j (\bb_i,\bb_j) \geq \sum_{i=1}^{k-1} \alpha_i^2 - 2c_n \sum_{1 \leq i < j \leq k-1} |\alpha_i \alpha_j|,$$
and so
$$|\cos \theta\ | = \frac{|(\bx,\bb_k)|}{\|\bx\|} \leq \frac{\sum_{i=1}^{k-1} |\alpha_i| |(\bb_i,\bb_k)|}{\|\bx\|} \leq \frac{ c_n \sum_{i=1}^{k-1} |\alpha_i|}{\sqrt{\sum_{i=1}^{k-1} \alpha_i^2 - 2c_n \sum_{1 \leq i < j \leq k-1} |\alpha_i \alpha_j|}}.$$
We want this quantity to be $\leq 1/2$, which is equivalent to saying that
\begin{equation}
\label{cn1}
4 c_n^2 \left( \sum_{i=1}^{k-1} |\alpha_i| \right)^2 \leq \sum_{i=1}^{k-1} \alpha_i^2 - 2c_n \sum_{1 \leq i < j \leq k-1} |\alpha_i \alpha_j|.
\end{equation}
Manipulating~\eqref{cn1}, we obtain
\begin{equation}
\label{cn2}
f(\alpha_1,\dots,\alpha_{k-1}) := \frac{\sum_{1 \leq i < j \leq k-1} |\alpha_i \alpha_j|}{\sum_{i=1}^{k-1} \alpha_i^2} \leq \frac{1-4 c_n^2}{8c_n^2 + 2c_n}.
\end{equation}
In other words, $|\cos \theta\ | \leq 1/2$ if and only if~\eqref{cn2} holds for all~$\alpha_1,\dots,\alpha_{k-1} \in \real$. Hence we want to maximize $f(\alpha_1,\dots,\alpha_{k-1})$ and prove that this maximum is no bigger than the right hand side of~\eqref{cn2}. We can assume without loss of generality that all $\alpha_i$ are nonnegative. For every $1 \leq i \leq k-1$,
$$f_i := \frac{\partial}{\partial\alpha_i} f(\alpha_1,\dots,\alpha_{k-1}) = \frac{ \left( \sum_{j=1}^{k-1} \alpha_j^2 \right) \left( \sum_{j \neq i} \alpha_j \right) - 2 \alpha_i \left( \sum_{1 \leq l < j \leq k-1} \alpha_l \alpha_j \right)}{\left( \sum_{j=1}^{k-1} \alpha_j^2 \right)^2}.$$
Then $\baa = (\alpha_1,\dots,\alpha_{k-1})$ is a critical point of $f$ if and only if $f_i(\baa) = 0$ for all $1 \leq i \leq k-1$, which is equivalent to
\begin{equation}
\label{fi1}
\alpha_i = \frac{1}{2} \left( \frac{\sum_{j=1}^{k-1} \alpha_j^2}{\sum_{1 \leq i < j \leq k-1} \alpha_i \alpha_j} \right) \left( \sum_{j=1, j \neq i}^{k-1} \alpha_j \right) = \frac{1}{2 f(\baa)} \sum_{j=1, j \neq i}^{k-1} \alpha_j.
\end{equation}
Summing~\eqref{fi1} over all~$i$, we obtain: 
$$\sum_{i=1}^{k-1} \alpha_i = \frac{1}{2 f(\baa)} \sum_{i=1}^{k-1} \sum_{j=1, j \neq i}^{k-1} \alpha_j = \frac{k-2}{2 f(\baa)} \sum_{i=1}^{k-1} \alpha_i,$$
which means that $\baa$ is a critical point of $f$ if and only if $\frac{k-2}{2 f(\baa)} = 1$, i.e., $f(\baa) = \frac{k-2}{2}$. Notice that this happens when $\alpha_1 = \dots = \alpha_{k-1} \neq 0$:
$$f(\alpha_1,\dots,\alpha_1) = \frac{\binom{k-1}{2} \alpha_1^2}{(k-1) \alpha_1^2} = \frac{k-2}{2},$$
i.e., $f$ is constant on the line $\{ \alpha_1 = \dots = \alpha_{k-1} \} \setminus \{ 0 \}$. Computing the Hessian matrix of $f$ at any point~$\baa$ with equal positive coordinates, we obtain:
$$H(f) = \frac{1}{k-1} \begin{pmatrix} -(k-2) & 1 & \dots & 1 \\ 1 & -(k-2) & \dots & 1 \\ \vdots & \vdots & \ddots & \vdots \\ 1 & 1 & \dots & -(k-2) \end{pmatrix}.$$
It is a $(k-1) \times (k-1)$ symmetric matrix with a simple eigenvalue $0$ and eigenvalue~$-1$ of multiplicity~$k-2$. Hence, as a symmetric bilinear form on the tangent space to the graph of~$f$ along the entire line $\{ \alpha_1 = \dots = \alpha_{k-1} \} \setminus \{ 0 \}$, $H(f)$ is negative semidefinite with $1$-dimensional kernel $\spn_{\real} \{ (1,\dots,1)^{\top} \}$. This implies that $f$ assumes its maximum at $(\alpha,\dots,\alpha) \neq \bo$. Hence we have $|\cos \theta\ | \leq 1/2$ if and only if
$$\frac{k-2}{2} \leq \frac{1-4 c_n^2}{8c_n^2 + 2c_n},$$
for all $k \leq n$. This is equivalent to saying that
$$4(n-1) c_n^2 + (n-2) c_n - 1 \leq 0.$$
For positive $c_n$, equality in this inequality holds if and only if $c_n$ is as in~\eqref{cn_def}.
\endproof

\begin{rem}
Notice that~$c_n$ defined in~\eqref{cn_def} is not much smaller than~$1/n$. For instance, for $n=1000$, $c_n = 0.00099801587...$ as compared to~$1/n = 0.001$. Furthermore,
\begin{eqnarray*}
\lim_{n \to \infty} \left( c_n/(1/n) \right) & = & \left( \lim_{n \to \infty} \frac{n}{8(n-1)} \right) \lim_{n \to \infty} \left( \sqrt{(n-2)^2 + 16(n-1)} - (n-2) \right) \\
& = & \frac{1}{8} \lim_{n \to \infty} \frac{16(n-1)}{\sqrt{(n-2)^2 + 16(n-1)} + (n-2)} = 1.
\end{eqnarray*}
This suggests that asymptotically as $n \to \infty$ the family of lattices~$A_n^*$ of Example~\ref{ETF} comes as close as possible to~$\W_n^*$ with respect to coherence. 
\end{rem}
\medskip

Approximation of WR lattices with respect to coherence has previously been considered in~\cite{WR_sim}, where a sequence of integer lattices approximating the hexagonal lattice in the plane was constructed (Theorem~1.6). Here we also construct an infinite family of integral WR planar lattices with arbitrarily small coherence and controlled minimal norm and denominator, thus approximating~$\zed^2$. This is a result in the spirit of Diophantine approximation.

\begin{prop} Let $0 < \eps \leq 1/2$ and $D$ be a positive squarefree integer. There exists an integral well-rounded lattice 
\begin{equation}
\label{L_def}
L = \begin{pmatrix} \sqrt{q} & p/\sqrt{q} \\ 0 & r\sqrt{D}/\sqrt{q} \end{pmatrix} \zed^2 \subset \real^2,
\end{equation}
where $p,q,r \in \zed_{>0}$ and $p^2+r^2D = q^2$, so that $0 < C(L) = p/q < \eps$ with 
$$q \leq  \frac{2D}{1-\eps} \left( \frac{1}{\eps} + 2 \sqrt{\frac{1}{\eps} - 1} \right).$$
For this $L$, we have
$$|L| = \sqrt{q} \leq \sqrt{ \frac{2D}{1-\eps} \left( \frac{1}{\eps} + 2 \sqrt{\frac{1}{\eps} - 1} \right)}$$
and
$$\det(L) = r\sqrt{D} \leq 2\sqrt{D} \left( \sqrt{ \frac{1}{\eps^2} - 1} + \sqrt{ \frac{1}{\eps} + 1} \right).$$
\end{prop}

\proof
Let $\gamma = m/(n\sqrt{D}) > 1$ be a rational multiple of $1/\sqrt{D}$ such that
\begin{equation}
\label{gm}
\gamma \leq \sqrt{\frac{1+ \eps}{1 - \eps}}.
\end{equation}
We can assume without loss of generality that $\gcd(m,n) = 1$. Then
$$\sqrt{D} < m/n \leq \sqrt{3D}.$$
Let $p = m^2 - Dn^2$, $q = m^2 + Dn^2$ and $r = 2mn$, then $p^2+r^2D = q^2$ and
$$\frac{p}{q} = \frac{\gamma^2 Dn^2-Dn^2}{\gamma^2Dn^2+Dn^2} = \frac{\gamma^2-1}{\gamma^2+1} \leq \eps.$$
Further,
\begin{equation}
\label{q_bnd}
q = m^2 + Dn^2 = (\gamma^2 + 1) Dn^2,
\end{equation}
and, by~\eqref{gm},
$$\frac{m}{n} \leq \frac{\sqrt{D(1+\eps)}}{\sqrt{1-\eps}} = \frac{\sqrt{D(\frac{1}{\eps} + 1)}}{\sqrt{\frac{1}{\eps} - 1}}.$$
We can then take $m = \left[ \sqrt{D(\frac{1}{\eps} + 1)} \right]$ and $n = \left[ \sqrt{\frac{1}{\eps} - 1} \right] + 1$. Combining this observation with~\eqref{q_bnd} and~\eqref{gm}, we obtain
\begin{eqnarray*}
q & \leq & \left( \frac{1+ \eps}{1 - \eps} + 1 \right) D \left( \left[ \sqrt{\frac{1}{\eps} - 1} \right] + 1 \right)^2 \\
& \leq & \frac{2D}{1-\eps} \left( \frac{1}{\eps} + 2 \sqrt{\frac{1}{\eps} - 1} \right).
\end{eqnarray*}
Then the lattice $L$ as in \eqref{L_def} is equal to the scalar factor $\frac{1}{\sqrt{q}}$ times the lattice $\Omega(p,q)$ in Proposition~1.1 of~\cite{fletcher-jones}, and thus its asserted properties follow immediately from this proposition.
\endproof
\bigskip

\section{Conclusions}
\label{conclude}

Our main contribution in this paper was a detailed study of the properties of well-rounded nearly orthogonal and weakly nearly orthogonal lattices, respectively the sets $\W^*_n$ and $\W_n$. We showed that these sets do not contain any local maxima of the packing density function, and no minima either except for the integer lattice $\zed^n$. Specifically, there are no perfect lattices in $\W_n$, although there are many eutactic and even strongly eutactic lattices. This observation leads, in particular, to a conclusion that eutactic and strongly eutactic lattices are not necessarily minima or maxima of the packing density function~$\delta$, even when restricted to the space of well-rounded lattices. We are not aware of this fact previously recorded in the literature. This conclusion should be compared to the result of A. Ash  \cite{ash1} asserting that critical points of~$\delta$ on the space of all lattices occur specifically at eutactic ones. 

While the starting point of our investigation is the concept of nearly orthogonal lattices and a basic property of the nearly orthogonal bases as established in~\cite{baraniuk}, we move far beyond the results of~\cite{baraniuk} in the study of the geometric and optimization properties of these lattices. In addition to the above-mentioned results, we give rather explicit information on their numbers of minimal vectors and construct explicit families of examples. Further, we show that lattices in $\W_n^*$ possess a rather rare property of being spanned by any collection of $n$ linearly independent minimal vectors. One of our main objects of study is coherence of these lattices, as defined by analogy with coherence of frames. In particular, we establish an asymptotically sharp bound on coherence of a basis in $\real^n$ that forces its corresponding lattice to be in $\W_n^*$. 

On the other hand, coherence $C(L)$ does not seem to fully capture the angular separation of the minimal set of vectors of $L$: indeed, coherence of a frame is often called the {\it worst case coherence} since it really measures just the largest absolute value of cosine. In many situations, however, it may be sufficient to have small coherence on average, not necessarily worst case. Indeed, {\it average coherence} on frames has been investigated by several authors (see \cite{bajwa}, \cite{mixon}). In our context, one can define average coherence of a lattice $L$ to be
$$\A(L) := \frac{1}{|S'(L)|-1} \max \left\{ \sum_{\bwy \in S'(L) \setminus \{\bx\}} \frac{ \left| \left( \bx, \bwy \right) \right| }{\|\bx\| \|\bwy\| } : \bx \in S'(L) \right\},$$
where $S'(L)$ is a subset of $S(L)$ consisting of one vector from each $\pm$ pair: the specific choice is not important due to the absolute value on the summation terms in our definition. One possible avenue for future research would be to investigate properties of average coherence and its correlation with other lattice functions (e.g. packing density) on a variety of geometrically significant classes of lattices, including our sets $\W_n$ and~$\W^*_n$.
\bigskip

{\bf Acknowledgement:} We thank the anonymous referees for many remarks and suggestions that improved the quality of exposition.

\bigskip

\bibliographystyle{plain}  

\end{document}